\documentclass[11pt,leqno]{article}
\usepackage{amssymb}
\usepackage{color}
\newtheorem{theorem}{Theorem}[section]

\newtheorem{lemma}{Lemma}[section]
\newtheorem{corollary}{Corollary}[section]
\newtheorem{remark}{Remark}[section]

\newtheorem{definition}{Definition}[section]
\newcommand{\dfrac}{\displaystyle\frac}
\newcommand{\dint}{\displaystyle \int}
\renewcommand{\H}{{\cal H}}

\newcommand{\ve}{{\varepsilon}}

\definecolor{darkred}{rgb}{0.8,0,0}

\begin{document}
\title{\bf Orbital Stability for the Schr\"{o}dinger Operator Involving Inverse Square Potential }
\author{\Large Georgios P. Trachanas, Athens\footnote{gtrach@gmail.com, gtrach@math.ntua.gr} \\[6pt]
\Large Nikolaos B. Zographopoulos, Athens
\footnote{nzograp@gmail.com, zographopoulosn@sse.gr}
}
\date{}
\maketitle
\pagestyle{myheadings} \thispagestyle{plain} \markboth{}{}
\maketitle
\begin{abstract}
In this paper we prove the existence of orbitally stable standing waves for the critical Schr\"{o}dinger operator, involving potential of the form $\left(\frac{N-2}{2}\right)^2|x|^{-2}$. The approach, being purely variational, is based on the precompactness of any minimizing sequence with respect to the associated energy. Moreover, we discuss the case of the presence of a Hardy energy term, in conjunction with the behavior of the standing waves at the singularity.
\end{abstract}
%
%
%
%
\section{Introduction}
\setcounter{equation}{0}
We are concerned with the concept of orbital stability of the standing wave solutions for the following semilinear Schr\"{o}dinger equation
\begin{equation}\label{original.problem}
i\psi_t+\Delta \psi+\left(\frac{N-2}{2}\right)^2\frac{\psi}{|x|^2}+|\psi|^{q-2}\psi=0,\;\;\;\;\;\;\; x \in  \mathbb{R}^N,
\end{equation}
where $\psi:\mathbb{R}^N\times \mathbb{R}\rightarrow \mathbb{C},\; N\geq 3$ and $2<q<\frac{4}{N}+2$.

A standing wave is a solution of the form
\[
\psi(x,t)=e^{i\lambda t}u(x),\;\;\; \lambda \in \mathbb{R},\;\; u:\mathbb{R}^N\rightarrow \mathbb{R},
\]
therefore, (\ref{original.problem}) is reduced to the stationary equation
\begin{equation}\label{original.stationary.eq}
-\Delta u-\left(\frac{N-2}{2}\right)^2\frac{u}{|x|^2}+\lambda u-|u|^{q-2}u=0,\;\;\;\;\;\;\; x \in  \mathbb{R}^N.
\end{equation}
It is well-known that Hardy's inequality
\[
\left(\frac{N-2}{2}\right)^2 \int_{\mathbb{R}^N}\frac{|u|^2}{|x|^2}dx \leq \int_{\mathbb{R}^N}|\nabla u|^2dx,\;\;\;\; \mbox{for all}\; u\in C_0^\infty(\mathbb{R}^N),\;\; N\geq 3,
\]
is closely related to elliptic and parabolic equations involving inverse square potentials. The optimal constant $c_{*}:=\left(\frac{N-2}{2}\right)^2$ is the natural borderline separating existence from nonexistence. We know from \cite{vz00} that the linear heat equation with an inverse square potential,
\begin{equation} \label{heat equation}
u_t=\Delta u+\frac{c}{|x|^2}u,
\end{equation}
admits a global solution for $c<c_{*}$, both in the bounded and the whole space case. On the other hand, there are no solutions for $c>c_{*}$, even locally in time, due to instantaneous blow-up. In the critical case $c=c_{*}$, a global solution is defined in \cite{vzog1,vzog2,vz00}, but the functional framework is more complicated.

There is a quite extent literature related to the Schr\"{o}dinger operator equipped with an inverse square potential, $H_c:=-\Delta-c|x|^{-2}$. Regarding the stationary problem, we shall refer to some papers related to the semilinear equation
\begin{equation}
H_cu=f(x,u).
\end{equation}
In \cite{dupaigne}, the author exhibits some existence, uniqueness and regularity results on various types of solutions for the case $0<c\leq c_{*}$, posed on a bounded domain containing the origin, where $f(x,u)=u^q+tg(x)$, $g\geq 0$ is a smooth, bounded function and $t>0$. The case $t=0$ on a ball of $\mathbb{R}^N$ appears in \cite{brezis.dupaigne}.

In the strictly subcritical case $c\in (0,c_{*})$, we refer indicatively to \cite{terracini} on some existence and nonexistence results in $D^{1,2}(\mathbb{R}^N)$, concerning the case $f(x,u)=u^{\frac{N+2}{N-2}}$. We also mention \cite{smets}, where the case $f(x,u)=K(x)u^{\frac{N+2}{N-2}}$ is considered, yielding a nonexistence result for $c\geq c_{*}$, and $K(x)$ is a positive and bounded weight. In \cite{deng.jin}, some existence results are obtained in $H^1(\mathbb{R}^N)$, where $f(x,u)=a(x)u+|u|^{\frac{4}{N-2}}u+g(x,u)$ for certain functions $a,\; g$. In the context of kernel estimates for the semigroup $e^{-tH_c},\;0<c<c_{*}$, we refer to \cite{davies}; see also \cite{bft} for the subcritical case of the potential. Analogue bounds on the Schr\"{o}dinger heat kernel for the case of the Laplace-Beltrami operator on some Riemmanian manifold can be found in \cite{zhang.only}. A result on the long-time behavior for the magnetic Schr\"{o}dinger operator with the critical potential is obtained in \cite{cazacu}. Finally, we refer to \cite{adimurthi.esteban} for the spectral analysis, in $D^{1,2}(\mathbb{R}^N)$, of the operator
\[
L_{q,\eta}u:=-\Delta u-c_{*}\frac{qu}{|x|^2}-\eta u,
\]
where $0\leq q\leq 1$ and $\eta$ satisfies certain growth conditions. We also mention \cite{ef} for some estimates on the moments of the negative eigenvalues of the Schr\"{o}dinger operator in the critical case.

In \cite{su.wang.willem}, the existence of radially symmetric ground state solutions is established for the following semilinear equation
\begin{equation}\label{plaplacian}
-\Delta u+V(|x|)u=Q(|x|)u^{q-1},\;\;\;\;\; x\in \mathbb{R}^N,
\end{equation}
where $V(|x|)$ is nonnegative and presents a singular behavior at the origin of order $|x|^{-s},s\in [2,N)$. Analogue results are obtained in \cite{su.wang.willem.jde}, by the same authors, for an equation similar to (\ref{plaplacian}), driven by the p-Laplacian. In both cases, the procedure is based on proving compact inclusions for certain weighted Sobolev spaces. Also, we refer to the variational approach of \cite{costa} for the case $V(|x|)\geq -(c_{*}-\alpha)|x|^{-2},\; \alpha>0,\; V(|x|)|x|^{2}\rightarrow +\infty$, whenever $|x|\rightarrow 0$ or $|x|\rightarrow \infty$. Finally, we refer to \cite{bellazzini.bonanno} for the case of a nonnegative potential $V(x)=l^2|x|^{-2}+|x|^{-a}$, where $a>0$ and $l\in \mathbb{Z}$.

In the time-dependent case, nonlinear equations of the general form
\begin{equation}\label{general.setting}
i\psi_t+\Delta \psi+V(x)\psi+f(x,\psi)=0,
\end{equation}
arise in various contexts of mathematical physics (see for example \cite{rose.weinstein} and the references therein). Perhaps, in many circumstances, the orbital stability of the standing wave solutions is the only extracted information on the asymptotic behavior of nonlinear phenomena governed by (\ref{general.setting}). When $V(x)\equiv 0$, we refer to the stability result in  \cite{cazenave.lions} (see also \cite{shibata}). In the same paper, a Hartree type equation is considered, equipped with a potential of the form $V(x)=\sum_1^m\frac{c_i}{|x-x_i|}$, where $c_1,...,c_m$ are positive constants, and $x_1,...,x_m$ are given poles. The case $V\in L^\infty(\mathbb{R}^N)$ can be found in \cite{bellazzini.visciglia,rose.weinstein}. We especially point out  \cite{genoud} where, among other, the existence and stability is obtained for a singular potential $V(x)\sim |x|^{-b},\;b\in (0,2)$, near zero. In \cite{debouard.fukuizumi,genoud.stuart,jeanjean}, the authors prove the orbital stability for weighted nonlinearities of certain singularity, that is, they consider the case $f(x,\psi)=Q(x)|\psi|^{q-1}\psi$, where $Q$ is allowed to behave like $|x|^{-b},b\in (0,2)$. Concerning singularities at infinity, we refer to the instability results in \cite{chen,chen.liu} for the case of harmonic potentials $V(x)=|x|^2$. A case of quasilinear equations is treated in \cite{guo.chen}. In the subcritical case $V(|x|)=c|x|^{-2},\; c<c_{*}$, we refer to \cite{zhang.zheng}, where the long-time behavior is studied by applying scattering theory in the $H^1(\mathbb{R}^N)$-setting. On some Strichartz type estimates for the linear Schr\"{o}dinger and the linear wave equation, respectively, equipped with a subcritical inverse square potential, we refer to \cite{burq}.

As it turns out, as far as we know, the results obtained here are new and they cover the critical inverse square potential, involving the best constant, on the whole space.

In Section \ref{sec.space.setting}, we formulate the functional framework. The proper norm for the setting of problems (\ref{original.problem}) and (\ref{original.stationary.eq}) is the following
\begin{equation}\label{right.norm.introd}
\left\|u\right\|^2_H:=\lim_{\varepsilon \downarrow 0}(I_{B_\varepsilon^c}(u)-\Lambda_{\varepsilon}(u)) + \left\|u\right\|^2_{L^2(\mathbb{R}^N)},
\end{equation}
where
\[
I_{B_\varepsilon^c}(u):=\int_{B_\varepsilon^c}|\nabla u|^2dx-\left(\frac{N-2}{2}\right)^2 \int_{B_\varepsilon^c}\frac{|u|^2}{|x|^2}dx,
\]
is the Hardy functional on the complement of the ball of radius $\varepsilon$ and centered at the origin, and the surface integral
\[
\Lambda_{\varepsilon}(u):=\frac{N-2}{2}\varepsilon^{-1}\int_{|x|=\varepsilon}|u|^2dS,
\]
represents the Hardy energy at the singularity. This unexpected energy formulation was first introduced in the paper \cite{vzog1} in order to overcome a functional difficulty emerged by the eigensolutions presenting the maximal singularity. Further applications to linear singular parabolic problems of the above consideration can be found in \cite{trach.zogr,vzog2}.

The present approach is based on the existence of critical points for the following functional
\[
E(u):=\frac{1}{2}\, \lim_{\varepsilon \downarrow 0}(I_{B_\varepsilon^c}(u)-\Lambda_{\varepsilon}(u))-\frac{1}{q}\int_{\mathbb{R}^N}|u|^qdx,
\]
on certain subsets of the energy space $H$ with prescribed $L^2$-norm. Namely, we define the $L^2$-sphere of $H$,
\[
\Gamma:=\left\{u\in H:\int_{\mathbb{R}^N}|u|^2dx=\gamma \right\},
\]
for a given $\gamma>0$, and we study the following minimization problem
\begin{equation}\label{original.minimization.problem}
\left \{
\begin{array}{ll}
u\in \Gamma,    \\
E(u)=\min \left\{E(z):z\in \Gamma,\; z\; \mbox{is radial}\right\},
\end{array}
\right.
\end{equation}
that is, we are interested for radial ground state solutions. The set of all global minimizers of (\ref{original.minimization.problem}) is denoted by
\[
S_\gamma:=\left\{u\in \Gamma:u\; \mbox{is a solution of}\; (\ref{original.minimization.problem}) \right\}.
\] \vspace{0.2cm}

In Section \ref{sec.radial.case}, we recall the notion of orbital stability and we prove the main result which is stated in the following:
\begin{theorem}\label{main.result}
Assume that $2<q<\frac{4}{N}+2,\;N\geq 3$, and $\gamma>0$ is a given constant. Then, the set $S_\gamma$ is orbitally stable.
\end{theorem}
Based on a proper transformation as in \cite{vzog1,vzog2}, we derive the behavior of the standing waves at the origin. More precisely, we establish that the minimizers of (\ref{original.minimization.problem}) behave exactly as $|x|^{-(N-2)/2}$ at the singularity, thus they do not belong to $H^1 (\mathbb{R}^N)$. In this case, we address the appearance of the Hardy singularity term at the origin.
\begin{theorem} \label{behavior.origin}
For each $\gamma>0$, every minimizer of (\ref{original.minimization.problem}) behaves at the origin like $|x|^{-(N-2)/2}$.
\end{theorem}
At this point we mention that the first nonexistence in $H_0^1 (\Omega)$ result was given in \cite{filippas.tertikas}, concerning the corresponding to (\ref{original.stationary.eq}) linear problem. The first exact description of the behavior at the singularity, and thus nonexistence in the classical Sobolev spaces, was given in \cite{vzog1,vzog2} for all minimizers. In almost all cases, this behavior corresponds to the appearance of the Hardy singularity term as a well-defined real number (positive or negative depending on the case); the Hardy functional is also well-defined and finite as a principal value (also positive or negative depending on the case). Exception to this is the case of the so called $k$-improved Hardy functional (see \cite{vzog2}), where both the associated Hardy functional and the Hardy singularity term tend to infinity.  Finally, we mention that results concerning the behavior at the singularity, by a different method, may be found in \cite{felli.ferrero}.

Furthermore, we confirm the orbital stability for the following Schr\"{o}dinger equation
\begin{equation}\label{equation.hardy.at.infty}
i|x|^{-4}w_t+\Delta w+\left(\frac{N-2}{2}\right)^2\, \frac{w}{|x|^2}=-h(x)\,w,\;\;\;\;\;\;\; x \in  \mathbb{R}^N,
\end{equation}
for some $h$ decaying to zero at infinity, 
which also presents a Hardy-type energy. The somehow unexpected fact is that this energy appears in an additive way to the total energy, as an effect that comes from infinity, and maybe represents the main part of it. The result is based on the arguments of \cite{vzog1, vzog2}. If we consider the corresponding to (\ref{equation.hardy.at.infty}) linear heat equation,
\begin{equation}\label{darkrn1}
\left\{\begin{array}{ccll}
|x|^{-4}\, u_t &=& \Delta u + \left(\frac{N-2}{2}\right)^2\, \dfrac{u }{|x|^2},\ &x \in \mathbb{R}^N,\; t>0,\\
u(x,0) &=& u_0(x),\;\; &\mbox{for}\;\; x \in \mathbb{R}^N, \nonumber \\
u(x,t)&\to& 0,\;\; &\mbox{for} \;\; |x| \to \infty,\;t>0\,,
\end{array}\right.
\end{equation}
then, using similarity variables, we may prove the existence of the Hardy-type energy. This energy comes from infinity, is additive to the total energy and constitutes the main part of it.

In Section \ref{sec.general.case}, we extend Theorem \ref{main.result} by removing the hypothesis of radial symmetry on $\psi$ and $u$. Nevertheless, in the general case, the cost we have to pay is the introduction of a certain weight function on the nonlinearity. Note that there is no possibility for a "non-weighted" Hardy-Sobolev inequality to hold, thus it is unclear to us if standard methods (e.g. see \cite{cazenave}) may be applied. On the other hand, in Subsection \ref{IHS} we prove a weighted (or improved) Hardy-Sobolev inequality.

More precisely, we prove the orbital stability of the following equation
\begin{equation}\label{general.original.problem}
i\psi_t+\Delta \psi+\left(\frac{N-2}{2}\right)^2\frac{\psi}{|x|^2}+g(x)|\psi|^{q-2}\psi=0,
\end{equation}
on $\mathbb{R}^N\times \mathbb{R}$, for weight functions $g:\mathbb{R}^N\rightarrow [0,+\infty)$ decaying, with a certain rate, to zero at infinity.
%
%
\section{Space Setting}\label{sec.space.setting}
\setcounter{equation}{0}
In this section, we introduce a more convenient function space in order to slide over the effect of the Hardy term $c_{*}|x|^{-2}$. We define the space $\H$ to be the completion of the $C^{\infty}_0 (\mathbb{R}^N)$-functions under the norm
\begin{equation}\label{normv}
||v||^2_{\H} = \int_{\mathbb{R}^N} |x|^{-(N-2)}\, |\nabla v|^2\, dx + \int_{\mathbb{R}^N} |x|^{-(N-2)}\, |v|^2\, dx.
\end{equation}
Notice that, setting
\begin{equation}\label{brezis.transf.}
u(x)=|x|^{-\frac{N-2}{2}}v(x),
\end{equation}
in (\ref{normv}), the first integral simplifies the Hardy functional, since
\[
\int_{\mathbb{R}^N} |x|^{-(N-2)}\, |\nabla v|^2\, dx=\int_{\mathbb{R}^N}|\nabla u|^2dx-\left(\frac{N-2}{2}\right)^2\int_{\mathbb{R}^N}\frac{|u|^2}{|x|^2}dx,
\]
at least for $C_0^\infty(\mathbb{R}^N)$-functions. In the sequel, we denote transformation (\ref{brezis.transf.}) by $u={\cal T}(v)$, which is an isometry between the spaces $X=L^2(\mathbb{R}^N)$ and $\widetilde{X}=L^2(\mathbb{R}^N,d\mu),\; d\mu=|x|^{-(N-2)}dx$. Our setting is based on this equivalence, as we shall see below.

It is clear that $\H$, equipped with the scalar product
\begin{equation}
(v_1,v_2)_{\cal H}:=\mbox{Re}\int_{\mathbb{R}^N}|x|^{-(N-2)}\nabla v_1 \nabla v_2dx+\mbox{Re}\int_{\mathbb{R}^N}|x|^{-(N-2)}v_1v_2dx,
\end{equation}
is a well-defined real Hilbert space. An equivalent definition of $\H$ is the following:
\begin{lemma}\label{density.first}
Every function $v$, such that $||v||_{\H} < \infty$, belongs to $\H$.
\end{lemma}
\emph{Proof.} We adapt the arguments of \cite{vzog1} in our case. First, we prove the result for $v \in L^{\infty}(\mathbb{R}^N)$. For $\varepsilon>0$ small enough and $t>0$, we define the cutoff function $\rho_{\varepsilon}(t) \in C_0 (\mathbb{R}_{+} \backslash \{0\})$, $0 \leq \rho_{\varepsilon} \leq 1$, as:
\[
\rho_{\varepsilon}(t)= \left \{
\begin{array}{lllll}
0,   &  0 < t < \varepsilon^2, \\ (\log 1/\varepsilon)^{-1}\, \log (t/\varepsilon^2),  & \varepsilon^2 < t < \varepsilon, \\
1,   & \varepsilon < t < \frac{1}{\varepsilon}, \\ (\log \varepsilon)^{-1}\, \log (t \varepsilon^2),   & \frac{1}{\varepsilon} < t < \frac{1}{\varepsilon^2},
\\ 0,   &  t > \frac{1}{\varepsilon^2}.
\end{array}
\right.
\]
For a fixed $v \in L^{\infty}(\mathbb{R}^N)$, $||v||_{\H} < \infty,$ we define
$v_{\varepsilon}(x) = \rho_\ve (|x|)\, v(x)$.  Note that
\[
|\nabla \rho_\ve(|x|)|=\frac{c_\ve}{|x|}, \qquad \qquad c_\ve:=(\log (1/\ve))^{-1},
\]
for $x \in A_{\varepsilon} \cup A_{\frac{1}{\varepsilon}}$, where $A_{\varepsilon} = \{x \in \mathbb{R}^N: \ \varepsilon^2 < |x| <\varepsilon \}$ and $A_{\frac{1}{\varepsilon}} = \{x \in \mathbb{R}^N: \ \varepsilon^{-1} < |x| < \varepsilon^{-2}  \}$, being zero otherwise.  Then, we have
\begin{eqnarray}\label{H1a}
||v_{\varepsilon} - v||^2_{\H} &\leq&
2\int_{A_\varepsilon \cup A_{\frac{1}{\varepsilon}}} |x|^{-(N-2)}\, |\nabla \rho_\ve (|x|)|^2\,
|v|^2\, dx \nonumber \\ && + 2 \int_{B_{\varepsilon} \cup B^c_{\frac{1}{\varepsilon}}} |x|^{-(N-2)}\,
(1-\rho_\ve)^2\, (|\nabla v|^2 + |v|^2)\, dx.
\end{eqnarray}
We will prove that letting $\varepsilon \downarrow 0$, the integrals in (\ref{H1a}) tend to zero. The only one that is
delicate is the first. We have
\[
\dint_{A_\varepsilon} |x|^{-(N-2)}\, |\nabla \rho_\ve (|x|)|^2\, |v|^2\, dx \le C\|v\|_\infty^2 \int_{\ve^2}^{\ve}c_\ve^2\frac{dr}{r}= C\|v\|_\infty^2(\log (1/\ve))^{-1},
\]
and
\[
\dint_{A_{\frac{1}{\varepsilon}}} |x|^{-(N-2)}\, |\nabla \rho_\ve (|x|)|^2\, |v|^2\, dx \le C\|v\|_\infty^2 \int_{\ve^{-1}}^{\ve^{-2}}c_\ve^2\frac{dr}{r}= C\|v\|_\infty^2(\log (1/\ve))^{-1},
\]
and these tend to zero as $\ve\to 0$.

Finally, we will prove that for any $v$ with $||v||_{\mathcal{H}} < \infty$, $v$ is approximated by a sequence of bounded functions in $\H$. Indeed, if we define $v_n$ as follows
\[
v_n(x) = \left \{
\begin{array}{ll}
v(x),   & \mbox{if} \;\;\; |v(x)|\leq n, \\
n,   & \mbox{if} \;\;\; |v(x)|>n,
\end{array}
\right.
\]
then, we have
\begin{equation}
||v_n||_{\H} = \int_{C_n} |x|^{-(N-2)}\, |\nabla v|^2\,dx + \int_{C_n} |x|^{-(N-2)}\, |v|^2\,dx< \infty.
\end{equation}
Now, it is clear that the sets $C_n=\{x\in \mathbb{R}^N: \ |v(x)|>n \}$
form a monotone family and the  measure tends to zero in the limit
$n\to \infty$. This means that the above integral goes to zero as
$n\to\infty$ and the proof is complete. $\blacksquare$

Now, the approximation of the functions $v_{\varepsilon}$, which vanish around the origin, by $C_0^{\infty} (\mathbb{R}^N \backslash \{0\})$ is standard.
\begin{corollary}\label{density.second}
The set $C_0^{\infty} (\mathbb{R}^N \backslash \{0\})$ is dense in $\H$.
\end{corollary}
We introduce the space $H$ to be the isometric space of $\mathcal{H}$ under $\mathcal{T}$, that is, $H$ is defined as the completion of the set
\[
\{u=|x|^{-\frac{N-2}{2}}v:\; v\in C_0^{\infty} (\mathbb{R}^N) \}=\mathcal{T}(C_0^{\infty} (\mathbb{R}^N)),
\]
under the norm
\begin{equation}\label{norm.for.u}
N^2(u):=\int_{\mathbb{R}^N}|x|^{-(N-2)}|\nabla(|x|^{\frac{N-2}{2}}u)|^2dx+\int_{\mathbb{R}^N}|u|^2dx.
\end{equation}
Then, the following is immediate:
\begin{corollary}
The set $C_0^\infty(\mathbb{R}^N \backslash \{0\})$ is dense in $H$.
\end{corollary}
Note also that
\begin{equation}
\left\|v\right\|_{{\cal H}}=\left\|u\right\|_H,\;\;\; \mbox{for all}\; u\in C_0^\infty(\mathbb{R}^N \backslash \{0\}).
\end{equation}
Arguing as in \cite{vzog1}, the norm of $H$ is related to the Hardy functional by means of the formula (\ref{right.norm.introd}).
%
\section{The Radial Case}\label{sec.radial.case}
\setcounter{equation}{0}
\subsection{Space Properties}
In the following, $2^{*}:=2N/(N-2)$ stands for the critical Sobolev exponent.
\begin{lemma} \label{compactradial}
Let $\H_r$ be the subspace of $\H$ consisting of radial functions. Then,
\[
\H_r \hookrightarrow L^{q}(\mathbb{R}^N,|x|^{-\frac{N-2}{2}q}dx),
\]
with compact inclusion for any $2 < q < 2^{*},\; N\geq 3$.
\end{lemma}
\emph{Proof.} Let $v_n (r),\; r=|x|$, be a bounded sequence of $C_0^{\infty} (\mathbb{R}^N \backslash \{0\})$- functions in $\H$; without loss of generality we assume that $v_n$ converges weakly to $0$. We will prove that
\begin{equation} \label{eq2.3.1}
\int_{\mathbb{R}^N} |x|^{-\frac{N-2}{2}q}\, |v_n|^q\, dx = \int_{B_1} |x|^{-\frac{N-2}{2}q}\, |v_n|^q\, dx + \int_{B_1^c} |x|^{-\frac{N-2}{2}q}\, |v_n|^q\, dx \to 0,
\end{equation}
as $n \to \infty$.
We claim that, for $v_n$ in $\H(B_1)$, it holds $v_n \in H^1(B_1 \subset \mathbb{R}^2)$. This is clear from the following
\[
||v_n||_{\H(B_1)} \sim \int_0^1 r\,  |v_n'|^2\, dr + \int_0^1 r\, |v_n|^2\, dr.
\]
Using now the compact imbedding
\begin{equation}\label{compimbH1}
H^1(B_1 \subset \mathbb{R}^2) \hookrightarrow L^p (B_1 \subset \mathbb{R}^2),\;\;\;\;\;\; \mbox{for any}\;\;\;\; p \in (1,+\infty),
\end{equation}
we obtain
\[
\int_{B_1} |x|^{-(N-2)}\, |v_n|^p\, dx \to 0,\;\;\;\;\;\; \mbox{for any}\;\;\;\; p \in (1,+\infty),
\]
as $n \to \infty$. Then, for $p$ large enough,
\begin{eqnarray*}
\int_{B_1} |x|^{-\frac{N-2}{2}q}\, |v_n|^q\, dx & = & \int_{B_1} |x|^{-(N-2)\frac{p-2}{2p}q}\, | |x|^{-\frac{N-2}{p}} v_n |^q\, dx \\
& \leq & \left ( \int_{B_1} |x|^{-(N-2)\frac{p-2}{2(p-q)}q}\, dx  \right )^{\frac{p-q}{p}}\, \left (  \int_{B_1} |x|^{-(N-2)}\, |v_n|^p\, dx   \right )^{\frac{q}{p}}.
\end{eqnarray*}
The first integral in the right hand side is finite for $2<q<2^{*}$. Therefore,
\begin{equation}
\int_{B_1}|x|^{-q\frac{N-2}{2}}|v_n|^qdx\rightarrow 0,\;\;\;\; \mbox{as}\; n\rightarrow \infty.
\end{equation}
In order to obtain the same for the second integral in (\ref{eq2.3.1}), we use the equivalence of the norm of $\H_r$ with the norm of $H^1_r (\mathbb{R}^2)$, the radial subspace of $H^1 (\mathbb{R}^2)$, and the compact imbedding
\[
H^1_r (\mathbb{R}^2) \hookrightarrow L^q (\mathbb{R}^2),
\]
for any $2 < q < 2^{*}$. Then,
\[
\int_{\mathbb{R}^N} |x|^{-(N-2)}\, |v_n|^q\, dx \to 0,
\]
and finally
\begin{equation} \label{eq2.3.1b}
\int_{B_1^c} |x|^{-\frac{N-2}{2}q}\, |v_n|^q\, dx \leq \int_{B_1^c} |x|^{-(N-2)}\, |v_n|^q\, dx \to 0,
\end{equation}
as $n \to \infty$. $\blacksquare$ \vspace{0.2cm}\\
By unitary equivalence, the result holds in terms of $u$, hence the inclusion
\[
H_r\hookrightarrow L^q(\mathbb{R}^N),\;\;\; 2<q<2^{*},
\]
is compact, where $H_r$ is the subspace of $H$ consisting of radial functions.

The following weighted interpolation inequality is immediate by setting $\alpha=\beta=\gamma=-\frac{N-2}{2}$ in the classical paper \cite{ckn}.
\begin{lemma}\label{lemma.ckn.radial}
Let $2<q<\frac{2N}{N-2}$ and $N\geq 3$. Then, there exists a constant $C=C(N,q)>0$, such that
\begin{equation}\label{CKN.without.weight}
\int_{\mathbb{R}^N}|x|^{-q\frac{N-2}{2}}|v|^qdx\leq C\left(\int_{\mathbb{R}^N}|x|^{-(N-2)}|\nabla v|^2dx\right)^{\frac{N(q-2)}{4}}\left(\int_{\mathbb{R}^N}|x|^{-(N-2)}|v|^2dx\right)^{\frac{2q-N(q-2)}{4}},
\end{equation}
for any $v\in\mathcal{H}$.
\end{lemma}
\begin{remark}
Notice that inequality (\ref{CKN.without.weight}) holds without any hypothesis on radial symmetry.
\end{remark}
%
%
%
%
\subsection{Global Solution}
We quote the basic points of the standard theory of semilinear Schr\"{o}dinger equations (cf. \cite{cazenave,cazenave.haraux})    ensuring the global well-posedness on the initial value problem
\begin{equation}\label{IVP.in.subsection}
\left \{
\begin{array}{ll}
i\psi_t+\Delta \psi+\left(\frac{N-2}{2}\right)^2\frac{\psi}{|x|^2}+|\psi|^{q-2}\psi=0,    \\
\psi(0)=\psi_0.
\end{array}
\right.
\end{equation}
The operator defined by
\begin{equation}
\left \{
\begin{array}{ll}
D(\widetilde{{\cal L}}):=\{\varphi \in {\cal H}_r:\widetilde{{\cal L}}\varphi \in \widetilde{X}\},    \\
\widetilde{{\cal L}}\varphi:=|x|^{N-2}\nabla \cdot(|x|^{-(N-2)}\nabla \varphi),
\end{array}
\right.
\end{equation}
defines a self-adjoint operator with $\widetilde{{\cal L}}\leq 0$. Therefore, by unitary equivalence, the operator defined by
\begin{equation}
\left \{
\begin{array}{ll}
D({\cal L}):=\left\{\psi \in H_r:{\cal L} \psi \in X\right\},    \\
{\cal L}\psi:=\Delta \psi+\left(\frac{N-2}{2}\right)^2\frac{\psi}{|x|^2},
\end{array}
\right.
\end{equation}
is also self-adjoint with ${\cal L}\leq 0$, where $\psi=|x|^{-(N-2)/2}\varphi$. Furthermore, $i{\cal L}$ defines a skew-adjoint operator and generates a group of isometries on $H_r$.
Adapting step by step the arguments of \cite[Theorem 3.3.9]{cazenave}, we may easily conclude the following:
\begin{theorem}\label{thm.local.existence}
Let $N\geq 3$ and $2<q<\frac{4}{N}+2$. Then, for all $\psi_0 \in H_r$, there exist $T=T(\psi_0)>0$ and unique solution
\[
\psi \in C([0,T),H_r)\cap C^1([0,T),H_r^{-1}),
\]
of (\ref{IVP.in.subsection}), where $H^{-1}_{r}$ is the dual space of $H_r$. In addition, for all $t\in [0,T)$, the following properties hold:
\begin{equation}\label{cons.law.charge}
\int_{\mathbb{R}^N}|\psi(t)|^2dx=\int_{\mathbb{R}^N}|\psi_0|^2dx \;\;\;\; \mbox{(conservation of charge)},
\end{equation}
and
\begin{equation}\label{cons.law.energy}
E(\psi(t))=E(\psi_0) \;\;\;\; \mbox{(conservation of energy)}.
\end{equation}
\end{theorem}
Set
\[
F(u):=\frac{1}{q}\int_{\mathbb{R}^N}|u|^qdx.
\]
Interpreting inequality (\ref{CKN.without.weight}) in terms of $u$, it follows that there exists $\varepsilon \in (0,1)$, such that
\begin{equation}
F(u)\leq \frac{1-\varepsilon}{2}\left\|u\right\|_{H}^2+C(\left\|u\right\|_{L^2(\mathbb{R}^N)}),\;\;\; \mbox{for all}\; u\in H_r.
\end{equation}
Therefore, in analogy to \cite[Theorem 3.4.1]{cazenave}, we may set $T(\psi_0)=\infty$ in Theorem \ref{thm.local.existence}, for all $\psi_0 \in H_r$.
%
%
%
%
%
\subsection{Stability}\label{stability}
We introduce the functional $J:{\cal H}_r\rightarrow \mathbb{R}$ defined by
\[
J(v):=E(|x|^{-\frac{N-2}{2}}v)+\frac{1}{2}\int_{\mathbb{R}^N}|x|^{-(N-2)}|v|^2dx.
\]
Clearly, the minimization problem (\ref{original.minimization.problem}) is equivalent to the following,
\begin{equation}\label{minimum.J.v.radial}
\left \{
\begin{array}{ll}
v\in \widetilde{\Gamma},    \\
J(v)=\min \{J(z):z\in \widetilde{\Gamma} \},
\end{array}
\right.
\end{equation}
where
\[
\widetilde{\Gamma}:=\{v\in {\cal H}_r:|x|^{-\frac{N-2}{2}}v\in \Gamma \}.
\]
Both minimization problems are well-defined since, by inequality (\ref{CKN.without.weight}) and the property $2<q<\frac{4}{N}+2$, there exist $\delta>0$ and $K<\infty$ such that
\begin{equation}
J(v)\geq \delta \left\|v\right\|^2_{{\cal H}_r}-K,\;\;\; \mbox{for all}\; v\in \widetilde{\Gamma}.
\end{equation}
Let
\[
\widetilde{k}_\gamma:=\inf \{J(v):v\in \widetilde{\Gamma} \}.
\]
In the following lemma we establish the precompactness of any minimizing sequence with respect to $\widetilde{k}_\gamma$.
\begin{lemma}\label{precompactness.property}
Let $\gamma>0$. Then, any sequence $(v_n)$ in ${\cal H}_r$ satisfying
\begin{equation}
J(v_n)\rightarrow \widetilde{k}_\gamma,\;\;J'(v_n)\rightarrow 0,\;\; v_n\in \widetilde{\Gamma},
\end{equation}
contains a convergent subsequence in ${\cal H}_r$. Furthermore, its limit solves the minimization problem (\ref{minimum.J.v.radial}).
\end{lemma}
\emph{Proof.} First, observe that for $n$ large enough,
\begin{eqnarray}
\widetilde{k}_\gamma+o(1) & \geq & J(v_n)-\frac{1}{q}<J'(v_n),v_n> \nonumber \\
& = & \frac{q-2}{2q}\left\|v_n\right\|_{{\cal H}_r}^2,
\end{eqnarray}
that is, $\left\|v_n\right\|_{{\cal H}_r}$ is bounded. Going if necessary to a subsequence, still denoted by $(v_n)$, there exists $v\in {\cal H}_r$ such that $v_n\rightharpoonup v$ in ${\cal H}_r$. By Lemma \ref{compactradial}, $v_n\rightarrow v$ in $L^q(\mathbb{R}^N,|x|^{-q\frac{N-2}{2}}dx)$.\\
Also, notice that
\begin{eqnarray}
\left\|v_n-v\right\|_{{\cal H}_r}^2 & = & <J'(v_n)-J'(v),v_n-v> \nonumber \\
&& +\: \int_{\mathbb{R}^N}|x|^{-q\frac{N-2}{2}}(|v_n|^{q-2}v_n-|v|^{q-2}v)(v_n-v)dx.
\end{eqnarray}
Evidently,
\begin{equation}
<J'(v_n)-J'(v),v_n-v>\; \rightarrow \; 0,\;\;\; \mbox{as}\; n\rightarrow \infty.
\end{equation}
Moreover, it follows by H\"{o}lder's inequality that
\begin{equation}
\left|\int_{\mathbb{R}^N}|x|^{-q\frac{N-2}{2}}(|v_n|^{q-2}v_n-|v|^{q-2}v)(v_n-v)dx\right|\leq C\left(\int_{\mathbb{R}^N}|x|^{-q\frac{N-2}{2}}|v_n-v|^qdx\right)^{1/q},
\end{equation}
which tends to zero, as $n \rightarrow \infty$. Therefore, we found a $v\in {\cal H}_r$ and a subsequence $(v_n)$ such that
\begin{equation}
v_n\rightarrow v,\;\;\; \mbox{in}\; {\cal H}_r.
\end{equation}
It is obvious that $v\in \widetilde{\Gamma}$. Finally, by the weak lower semicontinuity of $J$ and the definition of $\widetilde{k}_\gamma$, we obtain $J(v)=\widetilde{k}_\gamma$, and the proof is complete. $\blacksquare$
\vspace{0.1cm}

By a solution to (\ref{original.stationary.eq}) we mean a couple $(\lambda_\gamma,u_\gamma)\in \mathbb{R}\times H_r$, where $\lambda_\gamma$ is the Lagrange multiplier associated to the critical point $u_\gamma$ of $E$ on $\Gamma$. If $v_\gamma$ is a global minimizer of $J$ on $\widetilde{\Gamma}$, then there exist $(\lambda_\gamma,v_\gamma)\in \mathbb{R}\times {\cal H}_r$ solving the elliptic equation
\begin{equation}\label{elliptic.eq.in.vgamma}
-\nabla \cdot (|x|^{-(N-2)}\nabla v_\gamma)+\lambda_\gamma |x|^{-(N-2)}v_\gamma-|x|^{-q\frac{N-2}{2}}|v_\gamma|^{q-2}v_\gamma=0.
\end{equation}
By unitary equivalence, $S_\gamma \neq \emptyset$ and $e^{i\lambda_\gamma t}|x|^{-\frac{N-2}{2}}v_\gamma$ corresponds to a standing wave of (\ref{original.problem}). Observe that if $u_\gamma \in S_\gamma$, then $e^{i\lambda_\gamma t}u_\gamma \in S_\gamma$, for all $t\geq 0$. Also, $e^{i\lambda_\gamma t}u_\gamma$ is a periodic function in time, hence we may say that $S_\gamma$ consists of a class of closed orbits. In this context, the following definition of (local) orbital stability makes sense.
\begin{definition}
The set $S_\gamma$ is said to be orbitally stable, if for any $\varepsilon>0$, there exists a $\delta>0$, such that for any global solution $\psi(t)$ of (\ref{IVP.in.subsection}) with
\begin{equation}
dist(\psi_0,S_\gamma)<\delta,
\end{equation}
it holds that
\begin{equation}
dist(\psi(t),S_\gamma)<\varepsilon,\;\;\;\; \mbox{for all}\; t\geq 0,
\end{equation}
where
\begin{equation}
dist(w,S_\gamma):=\inf_{z\in S_\gamma}\left\|w-z\right\|_{H_r}.
\end{equation}
\end{definition}
\textit{Proof of Theorem \ref{main.result}.} Arguing by contradiction, assume that there exist sequences $(\psi_{0n})\subset H_r,\;(t_n)\subset \mathbb{R}_{+}$, and $\varepsilon_0>0$ with
\begin{equation}\label{contadiction.argument.a}
\left\|\psi_{0n}-u\right\|_{H_r} \rightarrow 0,
\end{equation}
for some $u\in S_\gamma$, and such that the global solutions $\psi_n$ with initial values $\psi_{0n}$ satisfy
\begin{equation}\label{contadiction.argument.b}
\inf_{z\in S_\gamma}\left\|\psi_n(t_n)-z\right\|_{H_r} \geq \varepsilon_0.
\end{equation}
Let us set
\[
u_n:=\psi_n(t_n).
\]
Clearly, if
\[
k_\gamma:=\min \left\{E(u):u\in \Gamma,\; u\; \mbox{is radial} \right\},
\]
then
\[
\widetilde{k}_\gamma=k_\gamma+\frac{\gamma}{2}.
\]
We deduce from (\ref{contadiction.argument.a}) that
\begin{equation}\label{convergence.of.initial.values}
E(\psi_{0n})\rightarrow k_\gamma \;\;\;\; \mbox{and}\;\;\;\; \int_{\mathbb{R}^N}|\psi_{0n}|^2dx\rightarrow \gamma.
\end{equation}
Applying the conservation laws (\ref{cons.law.charge}) and (\ref{cons.law.energy}), we obtain
\begin{equation}\label{convergence.of.maximal.solutions}
E(u_n)\rightarrow k_\gamma \;\;\;\; \mbox{and}\;\;\;\; \int_{\mathbb{R}^N}|u_n|^2dx\rightarrow \gamma.
\end{equation}
We choose a sequence $\beta_n\rightarrow 1$, such that
\begin{equation}
E(\beta_nu_n)\rightarrow k_\gamma \;\;\;\; \mbox{and}\;\;\;\; \int_{\mathbb{R}^N}|\beta_nu_n|^2dx=\gamma,
\end{equation}
that is, $\beta_nu_n$ is a minimizing sequence of problem (\ref{original.minimization.problem}). For example, we may insert
\begin{equation}
\beta_n=\left(\frac{\gamma}{\int_{\mathbb{R}^N}|u_n|^2dx}\right)^{1/2}.
\end{equation}
Setting $v_n=|x|^{\frac{N-2}{2}}u_n$, we obtain
\begin{equation}
J(\beta_nv_n)\rightarrow \widetilde{k}_\gamma \;\;\;\; \mbox{and}\;\;\;\; \int_{\mathbb{R}^N}|x|^{-(N-2)}|\beta_nv_n|^2dx=\gamma.
\end{equation}
Ekeland's variational principle \cite[Theorem 2.4]{willem} yields to the existence of another sequence $\zeta_n\in \widetilde{\Gamma}$ satisfying
\begin{equation}
J(\zeta_n)\rightarrow \widetilde{k}_\gamma,\;\;\; J'(\zeta_n)\rightarrow 0,\;\;\; \left\|\zeta_n-\beta_nv_n\right\|_{{\cal H}_r}<\frac{1}{n}.
\end{equation}
By Lemma \ref{precompactness.property}, we deduce that there exist a subsequence in ${\cal H}_r$, still denoted by $\zeta_n$, and a $\zeta \in S_\gamma$ such that $\zeta_n\rightarrow |x|^{\frac{N-2}{2}}\zeta$ in ${\cal H}_r$.
Therefore,
\begin{eqnarray}
\inf_{z\in S_\gamma}\left\|u_n-z\right\|_{H_r} & \leq & ||v_n-|x|^{\frac{N-2}{2}}\zeta ||_{{\cal H}_r} \nonumber \\
& = & ||v_n-\beta_nv_n+\beta_nv_n-\zeta_n+\zeta_n-|x|^{\frac{N-2}{2}}\zeta ||_{{\cal H}_r} \nonumber \\
& \leq & |1-\beta_n|\left\|v_n\right\|_{{\cal H}_r}+\frac{1}{n}+||\zeta_n-|x|^{\frac{N-2}{2}}\zeta ||_{{\cal H}_r}\rightarrow 0,
\end{eqnarray}
which contadicts (\ref{contadiction.argument.b}) and the proof is complete. $\blacksquare$
\vspace{0.1cm}
%
\subsection{Behavior at the Origin}
In this section, we describe the behavior of the standing wave solutions, obtained in Subsection \ref{stability}, around the origin. \vspace{0.2cm} \\
\emph{Proof of Theorem \ref{behavior.origin}}.\ \ Consider a global minimizer $u_\gamma \in S_\gamma$ of problem (\ref{original.minimization.problem}) for some $\gamma>0$, which can be chosen to be nonnegative on $\mathbb{R}^N \backslash \{0\}$. Then, there exists a Lagrange multiplier $\lambda_\gamma$, such that the pair $(\lambda_\gamma,v_\gamma)\in \mathbb{R}\times \mathcal{H}_r$ solves the stationary problem (\ref{elliptic.eq.in.vgamma}), where $v_\gamma=|x|^{(N-2)/2}u_\gamma$. \\
We introduce the following transformation
\begin{equation}
\widetilde{w}_\gamma(t)=v_\gamma(r),\;\;\;\;\; t=(-\log r)^{-\frac{1}{N-2}},\;\; r=|x|,
\end{equation}
for $0<r<1$. Note that
\[
\lim_{t\rightarrow +\infty}\widetilde{w}_\gamma(t)=\lim_{r\rightarrow 1^{-}}v_\gamma(r).
\]
Then, $\widetilde{w}_\gamma$ satisfies
\[
-\Delta \widetilde{w}_\gamma+\lambda (N-2)^2V_1(t)\widetilde{w}_\gamma+(N-2)^2V_2(t)|\widetilde{w}_\gamma|^{q-2}\widetilde{w}_\gamma=0,
\]
where
\[
V_1(t)=\exp \left(-2t^{-(N-2)}\right)t^{-2(N-1)} \;\;\; \mbox{and}\;\;\; V_2(t)=\exp \left(\left(q\frac{N-2}{2}-2\right)t^{-(N-2)}\right).
\]
If we set $V_1(0)=V_2(0)=0$, then $V_1,\; V_2$ are continuous functions. Standard regularity results imply that $\widetilde{w}_\gamma (0)$ is a well-defined real number, while Maximum Principle implies that this number cannot be zero. Thus $v_\gamma(0)$ is well-defined and positive, which means that $u_\gamma$ behaves at zero exactly as $|x|^{-(N-2)/2}$.\ $\blacksquare$
\vspace{0.1cm} \\
Consequently, the results of \cite[Section 2.5]{vzog1} yield the presence of a correcting term in the energy norm. The exact value of the norm of $u_\gamma$ is given by the formula
\[
\left\|u_\gamma\right\|^2_H=I_{\mathbb{R}^N}(u_\gamma)-\Lambda (u_\gamma)+\left\|u_\gamma\right\|^2_{L^2(\mathbb{R}^N)},
\]
where
\[\Lambda(u_\gamma)=\frac{N(N-2)}{2}\ \omega_N\ v^2_\gamma(0).
\]
On the contrary, we shall examine in next section a case where the new energy term acts in an additive way to the total energy.
\subsection{The Case of Hardy Energy at Infinity}
We consider the following Cauchy problem
\begin{equation}\label{IVP.hardy.at.infinity}
\left \{
\begin{array}{ll}
i|x|^{-4}w_t+\Delta w+\left(\frac{N-2}{2}\right)^2\frac{w}{|x|^2}=-|x|^{q(N-2)-2N}|w|^{q-2}w,    \\
w(0)=w_0,
\end{array}
\right.
\end{equation}
where $w=w(|x|),\; N\geq 3$ and $2<q<\frac{4}{N}+2$. The approach is based on the unitary equivalence with problem (\ref{IVP.in.subsection}). To this end, we introduce the Kelvin transformation in the form
\begin{equation}
\psi(y)=|x|^{N-2}w(x),\;\;\; x=\frac{y}{|y|^2},
\end{equation}
denoted by $\psi=\mathcal{K}(w)$. Clearly, for smooth functions, it holds that
\[
\Delta_y\psi(y)=|x|^{N+2}\Delta_xw(x) \;\;\; \mbox{and}\;\;\;\; \frac{\psi(y)}{|y|^2}=|x|^{N+2}\frac{w(x)}{|x|^2},
\]
and the equations are equivalent for $y\neq 0$. The differences appear in the energy near the singularity versus the energy at infinity.

As it was mentioned formerly, the energy space corresponding to (\ref{IVP.in.subsection}) is $H_r$ with norm given by (\ref{right.norm.introd}). We shall use this formulation for the definition of the energy space $\mathcal{W}$ corresponding to (\ref{IVP.hardy.at.infinity}). Indeed, following exactly the proposal stated in \cite{vzog1}, the space $\mathcal{W}$ is defined as the isometric space of $H_r$ under transformation $\mathcal{K}$ and the correct energy norm is given by the following formula
\[
\left\|w\right\|_{\mathcal{W}}^2=\lim_{\varepsilon \downarrow 0}(I_{1/\varepsilon}(w)+\Lambda_{1/\varepsilon}(w))+\int_{\mathbb{R}^N}|x|^{-4}|w|^2dx,
\]
where
\[
I_{1/\varepsilon}(w):=\int_{B_{1/\varepsilon}}|\nabla w|^2dx-\left(\frac{N-2}{2}\right)^2\int_{B_{1/\varepsilon}}\frac{|w|^2}{|x|^2}dx,
\]
and
\[
\Lambda_{1/\varepsilon}(w):=\frac{N-2}{2}\varepsilon \int_{|x|=\frac{1}{\varepsilon}}|w|^2dS.
\]
 More precisely, $\mathcal{W}$ is defined as the completion of the $C_0^\infty(\mathbb{R}^N \backslash \{0\})$-functions under the norm
\begin{equation}
\left\|w\right\|_{\mathcal{W}}^2=I_{\mathbb{R}^N}(w)+\int_{\mathbb{R}^N}|x|^{-4}|w|^2dx.
\end{equation}
\begin{remark}
The question of well-posedness, as well as, the existence and orbital stability of standing waves for (\ref{IVP.hardy.at.infinity}) in the space $\mathcal{W}$ is understood through the unitary equivalence with $H_r$.
\end{remark}
\begin{remark}
Note that in the case of a function $w$, which behaves at infinity like $|x|^{-(N-2)/2}$, we have again the appearance of a correction energy term with measure
\begin{equation}
\Lambda_{\infty}(w)=\frac{N(N-2)}{2}\ \omega_N\ v^2(0),
\end{equation}
where $v(x/|x|^2)=|x|^{\frac{N-2}{2}}w(x)$. However, there is a notable difference with problem (\ref{IVP.in.subsection}), since in this case the singularity effect acts in an additive way to the usual Hardy functional.
\end{remark}
%
\section{The General Case}\label{sec.general.case}
\setcounter{equation}{0}
In this section, we are concerned with the orbital stability of the Schr\"{o}dinger operator, equipped with the critical inverse square potential, without assuming any symmetry hypothesis. We consider the problem
\begin{equation}\label{IVP.in.nonradial.case}
\left \{
\begin{array}{ll}
i\psi_t+\Delta \psi+\left(\frac{N-2}{2}\right)^2\frac{\psi}{|x|^2}+g(x)|\psi|^{q-2}\psi=0,    \\
\psi(0)=\psi_0,
\end{array}
\right.
\end{equation}
where $g\in L^1(\mathbb{R}^N)$ is a nonnegative function. \\
The approach is an adaptation of the method used in the radial case. Therefore, we make an outline of the basic steps.
\begin{lemma} \label{compactv}
Let $1 \leq q < 2^{*}$. If $g$ satisfies
\begin{equation}\label{condition.on.weight.g}
g\sim r^\omega,\;\; \mbox{with}\; \left \{
\begin{array}{ll}
\omega>-N+\frac{q(N-2)}{2},   & \mbox{at}\;\;\; 0,   \\
\omega<-N+\frac{q(N-2)}{2},   & \mbox{at}\;\;\; \infty,
\end{array}
\right.
\end{equation}
then, the inclusion $\H \hookrightarrow L_g^q (\mathbb{R}^N,|x|^{-\frac{N-2}{2}q}dx)$ is compact, where
\[
L_g^q (\mathbb{R}^N,|x|^{-\frac{N-2}{2}q}dx):=\left\{v\in L^1(\mathbb{R}^N):\int_{\mathbb{R}^N}g|x|^{-q\frac{N-2}{2}}|v|^qdx<\infty \right\}.
\]
\end{lemma}
\emph{Proof.} Let $v_n$ be a bounded sequence of $C_0^{\infty} (\mathbb{R}^N \backslash \{0\})$- functions in $\H$; without loss of generality, we assume that $v_n$ converges weakly to $0$. Then, we split the norm of $v_n$ as
\[
||v_n||_{\H} = ||v_n||_{\H(B_1)} + ||v_n||_{\H(B_1^c)}.
\]
In addition, we study separately the radial parts $v_n^r=v_n^r(r),\; r=|x|$, and the nonradial parts $v_n^{nr}$ of $v_n$. The proof consists of four steps.

\emph{Step 1: The radial part $v_n^r$ in $B_1$}.
Arguing as in lemma \ref{compactradial}, we obtain that $v_n^r \in H^1(B_1 \subset \mathbb{R}^2)$, whenever  $v_n^r\in \H(B_1)$.    Therefore,
\[
\int_{B_1} |x|^{-(N-2)}\, |v_n^r|^p\, dx \to 0,\;\;\;\;\;\; \mbox{for any}\;\;\;\; p \in (1,+\infty),
\]
as $n \to \infty$. Then, for $p$ large enough,
\begin{eqnarray*}
\int_{B_1} |x|^{-\frac{N-2}{2}q}\, g\, |v_n^r|^q\, dx & = & \int_{B_1} |x|^{-(N-2)\frac{p-2}{2p}q}\, g\,  ||x|^{-\frac{N-2}{p}} v_n^r|^q\, dx \\
& \leq & \left ( \int_{B_1} |x|^{-(N-2)\frac{p-2}{2(p-q)}q}\, g^{\frac{p}{p-q}}\, dx  \right )^{\frac{p-q}{p}}\, \left (  \int_{B_1} |x|^{-(N-2)}\, |v_n^r|^p\, dx   \right )^{\frac{q}{p}}.
\end{eqnarray*}
Let $g \sim r^{\omega}$ at $0$. Then, the first integral of the right hand side is finite if
\begin{equation}\label{omegaB1r}
    \omega > (N-2)\frac{q}{2} - N,
\end{equation}
Finally, we get that
\[
\int_{B_1} |x|^{-\frac{N-2}{2}q}\, g\, |v_n^r|^q\, dx \to 0,\;\;\;\;\;\; \mbox{as}\;\;\;\; n \to \infty,
\]
for $g \sim r^{\omega}$, at $0$, and $\omega$ satisfying (\ref{omegaB1r}). \vspace{0.2cm} \\
\emph{Step 2: The nonradial part $v_n^{nr}$ in $B_1$}. For the nonradial parts $v_n^{nr}$ of $v_n$, we observe that
\begin{equation} \label{nonradialH1}
|x|^{-(N-2)/2}\, v_n^{nr} \in H^1 (\mathbb{R}^N).
\end{equation}
This follows as in \cite{vz00} or as in \cite[pg.196]{filippas.tertikas}; setting $u_n^{nr}=|x|^{-(N-2)/2}\, v_n^{nr}$, we have
\[
||v_n^{nr}||_{\H}^2 = I_{\mathbb{R}^N}(u_n^{nr}) + \int_{\mathbb{R}^N} |u_n^{nr}|^2\, dx.
\]
However, using decomposition into spherical harmonics,
\[
I_{\mathbb{R}^N}(u_n^{nr}) \geq c\, \int_{\mathbb{R}^N} |\nabla u_n^{nr}|^2\, dx,
\]
and (\ref{nonradialH1}) holds. Moreover,
\[
||v_n^{nr}||_{\H} \geq c\, ||v_n^{nr}||_{H^1(\mathbb{R}^N)},
\]
so $v_n^{nr}$ is also bounded in $H^1(\mathbb{R}^N)$.
Thus, $|x|^{-(N-2)/2}\, v_n^{nr} \in H^1 (B_1)$, with $H^1 (B_1)$ to be compactly embedded into $L^{p} (B_1)$, $1 \leq p < 2^{*}$.

For $\varepsilon >0$ small enough, we set
\begin{equation} \label{Aqe}
A_{q,\varepsilon} = \left ( 1 - \frac{q}{2^{*} - \varepsilon} \right )^{-1}.
\end{equation}
It is clear that $A_{q,\varepsilon}$ is increasing as a function of $\varepsilon$, so
\[
A_{q,\varepsilon} > \left ( 1 - \frac{q}{2^{*}} \right )^{-1},
\]
for any $\varepsilon >0$ small enough. Then,
\begin{eqnarray*}
\int_{B_1} |x|^{-\frac{N-2}{2}q}\, g\, |v_n^{nr}|^q\, dx &=& \int_{B_1} g\, ||x|^{-\frac{N-2}{2}} v_n^{nr}|^q\, dx \\
&\leq& \left ( \int_{B_1} g^{A_{q,\varepsilon}}\, dx  \right )^{\frac{1}{A_{q,\varepsilon}}}\, \left (  \int_{B_1} ||x|^{-\frac{N-2}{2}}\, v_n^{nr}|^{2^{*} -\varepsilon}\, dx   \right )^{\frac{q}{2^{*} -\varepsilon}}.
\end{eqnarray*}
Let $g \sim r^{\omega}$ at $0$. Then, the first integral of the right hand side is finite if (\ref{omegaB1r}) holds.
Finally, we get that
\[
\int_{B_1} |x|^{-\frac{N-2}{2}q}\, g\, |v_n^{nr}|^q\, dx \to 0,\;\;\;\;\;\; \mbox{as}\;\;\;\; n \to \infty,
\]
for $g \sim r^{\omega}$, at $0$, and $\omega$ satisfying (\ref{omegaB1r}). \vspace{0.2cm} \\
\emph{Step 3: The radial part $v_n^{r}$ in $B^c_1$}. In the case of the exterior domain $B_1^c$, both in the radial and the nonradial case, we use the following Kelvin transform. We set
\begin{equation}\label{kelvintrans}
v(x) = w(y),\;\;\;\;\;\; y=\frac{x}{|x|^2}.
\end{equation}
The determinant of the Jacobian of the Kelvin transformation in dimension $d \geq 2$ is equal to $-|x|^{2d}$ and
\[
\left | \nabla_x w \left ( \frac{x}{|x|^2} \right ) \right |^2 = |x|^{-4}\, |\nabla_y w (y)|^2.
\]
Then,
\begin{equation}\label{normB1w}
||v_n||^2_{\H(B_1^c)} = \int_{B_1} |y|^{-(N-2)}\, |\nabla w_n|^2\, dy + \int_{B_1} |y|^{-N-2}\, |w_n|^2\, dy.
\end{equation}
We restrict now in the radial case. For the radial functions $w_n^r$, we have
\[
\int_{B_1} |y|^{-(N-2)}\, |\nabla w_n^r|^2\, dy + \int_{B_1} |y|^{-N-2}\, |w_n^r|^2\, dy \geq ||w_n^r||^2_{H^1(B_1 \subset \mathbb{R}^2)}.
\]
Hence, as in step 1,
\[
|y|^{-\frac{N-2}{p}}\, w_n^r\;\;\;\;\; \mbox{converges to 0 in}\;\; L^p(B_1),\;\; p \geq 1.
\]
Moreover,
\begin{eqnarray*}
\int_{B_1^c} |x|^{-\frac{N-2}{2}q}\, g\, |v_n^{r}|^q\, dx & = & \int_{B_1} |y|^{\frac{N-2}{2}q -2N}\, g\, |w_n^r|^q\, dy  \\
& = & \int_{B_1} |y|^{\frac{N-2}{2}q -2N + \frac{N-2}{p}q}\, g\, ||y|^{-\frac{N-2}{p}}\, w_n^r|^q\, dy  \\
& \leq & \left ( \int_{B_1} |y|^{(N-2)\frac{p+2}{2(p-q)}q - \frac{2Np}{p-q}}\, g^{\frac{p}{p-q}}\, dx  \right )^{\frac{p-q}{p}}\,  \\
&& \left (  \int_{B_1} |x|^{-(N-2)}\, |v_n^r|^p\, dx   \right )^{\frac{q}{p}}.
\end{eqnarray*}
Let $g \sim |x|^{\omega}$ at $\infty$. Then, the first integral of the right hand side is finite if
\[
 - \omega > -q(N-2)\frac{p+2}{2p} + N + \frac{q}{p}N,
\]
or, for $p$ large enough,
\begin{equation}\label{omegaB1cr}
\omega < \frac{q(N-2)}{2} - N.
\end{equation}
Finally, we conclude
\[
\int_{B_1} |x|^{-\frac{N-2}{2}q}\, g\, |v_n^r|^q\, dx \to 0,\;\;\;\;\;\; \mbox{as}\;\;\;\; n \to \infty,
\]
for $g \sim r^{\omega}$, at $0$, and $\omega$ satisfying (\ref{omegaB1cr}). \vspace{0.2cm} \\
\emph{Step 4: The nonradial part $v_n^{nr}$ in $B^c_1$}. Using transformation (\ref{kelvintrans}), we have that
\begin{equation}\label{normRNw}
||v_n||^2_{\H} = \int_{\mathbb{R}^N} |y|^{-(N-2)}\, |\nabla w_n|^2\, dy + \int_{\mathbb{R}^N} |y|^{-N-2}\, |w_n|^2\, dy.
\end{equation}
Moreover, as in Step 2, we may obtain
\[
||v_n^{nr}||^2_{\H} \geq \int_{\mathbb{R}^N} |\nabla ( |y|^{-\frac{N-2}{2}}\, w_n^{nr} ) |^2\, dy + \int_{\mathbb{R}^N} |y|^{-4}\, ||y|^{-\frac{N-2}{2}}\, w_n^{nr}|^2\, dy.
\]
Restricting ourselves in $B_1$, we have
\[
|y|^{-\frac{N-2}{2}}\, w_n^{nr}\;\;\;\;\;\; \mbox{is bounded in}\;\;\; H^1(B_1),
\]
so
\[
|y|^{-\frac{N-2}{2}}\, w_n^{nr}\;\;\;\;\;\; \mbox{converges to 0 in}\;\;\; L^p(B_1),\;\; 1 \leq p < 2^{*}.
\]
Then,
\begin{eqnarray*}
\int_{B_1^c} |x|^{-\frac{N-2}{2}q}\, g\, |v_n^{nr}|^q\, dx & = & \int_{B_1} |y|^{\frac{N-2}{2}q -2N}\, g\, |w_n^{nr}|^q\, dy \\
& = &\int_{B_1} |y|^{(N-2)q -2N}\, g\, | |y|^{-\frac{N-2}{2}}\, w_n^{nr} |^q\, dy \\
& \leq & \left ( \int_{B_1}  ( |y|^{(N-2)q -2N}\, g )^{A_{q,\varepsilon}}\, dy  \right )^{\frac{1}{A_{q,\varepsilon}}}\, \\
&& \left (  \int_{B_1} | |y|^{-\frac{N-2}{2}}\, w_n^{nr} |^{2^{*} -\varepsilon}\, dy   \right )^{\frac{q}{2^{*} -\varepsilon}},
\end{eqnarray*}
where $A_{q,\varepsilon}$ is defined in (\ref{Aqe}). Let $g \sim |x|^{\omega}$ at $\infty$. Then, the first integral of the right hand side is finite if
\[
 - \omega > 2N - (N-2)q - \frac{N}{A_{q,\varepsilon}} > N - \frac{N-2}{2}q,
\]
or $\omega$ satisfies (\ref{omegaB1cr}). Finally, we get that
\[
\int_{B_1^c} |x|^{-\frac{N-2}{2}q}\, g\, |v_n^{nr}|^q\, dx \to 0,\;\;\;\;\;\; \mbox{as}\;\;\;\; n \to \infty,
\]
for $g \sim r^{\omega}$, at $\infty$, and $\omega$ satisfying (\ref{omegaB1cr}). $\blacksquare$
\begin{remark} \label{remcomp}
(i) Since
\begin{equation}\label{negative}
-N + \frac{N-2}{2}q < 0,
\end{equation}
for every $1 \leq q < 2^{*}$, we conclude that $g$ maybe constant at the origin. On the other hand, it should decay to zero at infinity. As we saw in Lemma \ref{compactradial}, this is not the case if we restrict on radial functions.

(ii) Consider the functions $u = |x|^{-(N-2)/2}\, v$. Then,
\[
\int_{\mathbb{R}^N} |x|^{-\frac{N-2}{2}q}\, g\, |v|^q\, dx = \int_{\mathbb{R}^N} g\, |u|^q\, dx.
\]
In the special case $q=2$, $g$ should be a function such that
\[
g \sim r^{\omega},\;\;\; \mbox{with} \; \left \{
\begin{array}{ll}
\omega > -2,   &  \mbox{at}\;\;\; 0, \\
\omega < -2,   &  \mbox{at}\;\;\; \infty.
\end{array}
\right.
\]
as it is natural from the Hardy's inequality.

(iii) Any function $g$ with $g \in L^{1} (\mathbb{R}^N) \cap L^{\frac{2^{*}}{2^{*} -q}} (\mathbb{R}^N)$ satisfies (\ref{condition.on.weight.g}).
\end{remark}
By Lemma \ref{lemma.ckn.radial} and the modified behavior of $g$ at zero, the following is immediate:
\begin{corollary}\label{ckn.nonradial}
Let $2< q<2^{*},\;N\geq 3$. If $g$ satisfies
\begin{equation}\label{condition.on.weight.g.modified}
g\sim r^\omega,\;\; \mbox{with}\; \left \{
\begin{array}{ll}
\omega \geq 0,   & \mbox{at}\;\;\; 0,   \\
\omega<-N+\frac{q(N-2)}{2},   & \mbox{at}\;\;\; \infty,
\end{array}
\right.
\end{equation}
then, there exists a constant $C=C(N,q)>0$, such that
\begin{eqnarray}
\int_{\mathbb{R}^N}g|x|^{-q\frac{N-2}{2}}|v|^qdx & \leq & \left(\int_{\mathbb{R}^N}|x|^{-(N-2)}|\nabla v|^2dx\right)^{\frac{N(q-2)}{4}}
\nonumber \\
&&
\left(\int_{\mathbb{R}^N}|x|^{-(N-2)}|v|^2dx\right)^{\frac{2q-N(q-2)}{4}},
\end{eqnarray}
for all $v\in {\cal H}$.
\end{corollary}
The energy functional naturally associated to (\ref{IVP.in.nonradial.case}) is defined by
\[
E_g(u):=\frac{1}{2}\lim_{\varepsilon \downarrow 0}(I_{B_\varepsilon^c}(u)-\Lambda_{\varepsilon}(u))-\frac{1}{q}\int_{\mathbb{R}^N}g|u|^qdx.
\]
\begin{theorem}
Let $N\geq 3$ and $2<q<\frac{4}{N}+2$. If $g$ satisies (\ref{condition.on.weight.g.modified}), then, for all $\psi_0 \in H$, there exists a unique solution
\[
\psi \in C([0,\infty),H)\cap C^1([0,\infty),H^{-1}),
\]
of (\ref{IVP.in.nonradial.case}), where $H^{-1}$ is the dual space of $H$. In addition, for all $t\geq 0$, the following conservation laws hold:
\begin{equation}
\int_{\mathbb{R}^N}|\psi(t)|^2dx=\int_{\mathbb{R}^N}|\psi_0|^2dx,
\end{equation}
and
\begin{equation}
E_g(\psi(t))=E_g(\psi_0).
\end{equation}
\end{theorem}
Consider the minimization problem
\begin{equation}\label{minimum.Eg.nonradial}
\left \{
\begin{array}{ll}
u\in \Gamma,    \\
E_g(u)=\min \left\{E_g(z):z\in \Gamma \right\},
\end{array}
\right.
\end{equation}
which is well-defined due to Corollary \ref{ckn.nonradial}. By using similar arguments to those used in Lemma \ref{precompactness.property}, we may obtain the existence of a standing wave solution $e^{i\lambda_\gamma t}u_\gamma \in H$ of (\ref{IVP.in.nonradial.case}), for a given $\gamma>0$.
If
\[
S_{g,\gamma}:=\left\{u\in H:u\; \mbox{is a solution of}\; (\ref{minimum.Eg.nonradial}) \right\}
\]
denotes the set of global minimizers, then a nonradial version of Theorem \ref{main.result} is stated in the following:
\begin{theorem}
Assume that $N\geq 3,\; 2<q<\frac{4}{N}+2$ and $g$ satisfies condition (\ref{condition.on.weight.g.modified}). Then, the set $S_{g,\gamma}$ is orbitally stable for problem (\ref{IVP.in.nonradial.case}).
\end{theorem}
\subsection{Improved Hardy-Sobolev Inequality} \label{IHS}
As already said in the Introduction, there is no possibility for a Hardy-Sobolev inequality to hold. This is clear since functions behaving at the origin like $|x|^{-(N-2)/2}$, do not belong to $L^{2^{*}}$. As in the bounded domain case, we have to consider Improved Hardy-Sobolev (IHS) inequalities; i.e. we have to find a weight function $h(x)$ such that
\begin{equation} \label{IHSRN}
||\phi||_H \geq \left ( \int_{\mathbb{R}^N} h(x)\, |\phi|^{2^{*}}  \right )^{(N-2)/N},\;\;\;\;\;\; \mbox{for any}\;\;\; \phi \in C_0^{\infty} (\mathbb{R}^N).
\end{equation}
Unfortunately, in the case of the whole space we cannot apply the method described in \cite[Lemma 9.1]{vzog2}. This method provides a way of constructing
optimal IHS inequalities (i.e. it provides $h(x)$) and it was used in \cite{trach.zogr, vzog2, zogr.jfa}. However, we are in position to present an IHS inequality, although it is not optimal. Before this, we remark the following:
\begin{enumerate}
  \item The nonradial part of $\phi$ belongs to $H^1 (\mathbb{R}^N)$ (see (\ref{nonradialH1})), so in this case (\ref{IHSRN}) holds with $h(x)=1$.
  \item The radial part $\phi_r$ of $\phi$ satisfies (\ref{IHSRN}) with $h(r)=|x|^2$. This may be obtained by following the arguments of Lemma \ref{compactradial}; in this case $|x|^{(N-2)/2} \phi_r \in H^1 (\mathbb{R}^2)$.
  \item For any $R>0$ and $\phi \in C_0^{\infty} (B_R)$, inequality (\ref{IHSRN}) holds with
\[
h(x) = \left (-\log \left( \frac{|x|}{D}\right) \right )^{-\frac{2(N-1)}{N-2}},\;\;\;\; D>R.
\]
This optimal (IHS) inequality was proved in \cite{filippas.tertikas}. In the case of radial functions, $D$ is allowed to be equal to $R$ (see \cite{zogr.jfa}). Note that, since $D$ or $R$ appear in the inequalities, these inequalities depend on the domain.
\end{enumerate}

Directly from 1. and 2. we have the following.
\begin{lemma}
Assume that $h$ is defined as
\[
h(x) = \left \{
\begin{array}{ll}
|x|^2,   &  |x|<1, \\
1,   &  |x| \geq 1.
\end{array}
\right.
\]
Then, inequality (\ref{IHSRN}) is true.
\end{lemma}
However, from the point of view of 3. this inequality is not optimal.

\vspace{1cm}

\noindent \textsc{Acknowledgment.} The first author was supported by $E\Lambda KE$, National Technical University of Athens (grant number 65193600).

\

%

{\small
\bibliographystyle{amsplain}

\begin{thebibliography}{10}
%
%
%
%
\bibitem{adimurthi.esteban} Adimurthi, M.J. Esteban, \emph{An improved Hardy-Sobolev inequality in $W^{1,p}$ and its application to Schr\"{o}dinger operators}, Nonlinear Differential Equations Appl. \textbf{12} (2005), 243-263.
%
\bibitem{bft} G. Barbatis, S. Filippas, A. Tertikas, \emph{Critical heat kernel estimates for Schr\"{o}}dinger operators via Hardy-Sobolev inequalities, J. Funct. Anal.  \textbf{208} (2004), 1-30.
%
\bibitem{bellazzini.bonanno} J. Bellazzini, C. Bonanno, \emph{Nonlinear Schr\"{o}dinger equations with strongly singular potentials}, Proc. Roy. Soc. Edinburgh Sect. A \textbf{140} (2010), 707-721.
%
\bibitem{bellazzini.visciglia} J. Bellazzini, N. Visciglia, \emph{On the orbital stability for a class of nonautonomous NLS}, Indiana Univ. Math. J. \textbf{59} (2010), 1211-1230.
%
%
\bibitem{brezis.dupaigne} H. Brezis, L. Dupaigne, A. Tesei, \emph{On a semilinear elliptic equation with inverse-square potential}, Selecta Math. (N.S.) \textbf{11} (2005), 1-7.
%
\bibitem{burq} N. Burq, F. Planchon, J. G. Stalker, A. S. Tahvildar-Zadeh, \emph{Strichartz estimates for the wave and Schr\"{o}dinger equations with the inverse-square potentail}, J. Funct. Anal. \textbf{203} (2003), 519-549.
%
\bibitem{ckn} L. Caffarelli, R. Kohn, L. Nirenberg, \emph{First
order interpolation inequalities with weights}, Compositio Math. \textbf{53} (1984), 259-275.
%
\bibitem{cazacu} C. Cazacu, D. Krejcirik, \emph{The Hardy inequality and the heat equation with magnetic field in any dimension},  arXiv:1409.6433v1.
%
\bibitem{cazenave} T. Cazenave, \emph{Semilinear Schr\"{o}dinger equations}, Courant Lecture Notes in Mathematics \textbf{10}, 2003.
%
\bibitem{cazenave.haraux} T. Cazenave, A. Haraux, \emph{An introduction to semilinear evolution equations}, Oxford Lecture Series in Mathematics and its Applications \textbf{13}, 1998.
%
\bibitem{cazenave.lions} T. Cazenave, P. L. Lions, \emph{Orbital stability of standing waves for some nonlinear Schr\"{o}dinger equations}, Comm. Math. Phys. \textbf{85} (1982), 549-561.
%
\bibitem{chen} J. Chen, \emph{On the inhomogeneous nonlinear Schr\"{o}dinger equation with harmonic potential and unbounded coefficient}, Czechoslovak. Math. J. \textbf{60} (2010), 715-736.
%
\bibitem{chen.liu} J. Chen, Y. Liu, \emph{Instability of standing waves to the inhomogeneous nonlinear Schr\"{o}dinger equation with harmonic potential}, Illinois J. Math. \textbf{52} (2008), 1259-1276.
%
%
\bibitem{costa} D. G. Costa, J. Marcos do \'O, K. Tintarev, \emph{Compactness properties of critical nonlinearities and nonlinear Schr\"{o}dinger equations}, Proc. Edinb. Math. Soc. \textbf{56} (2013), 427-441.
%
%
\bibitem{davies} E. B. Davies, B. Simon, \emph{$L^p$ norms of noncritical Schr\"{o}dinger semigroups},  J.Funct.Anal. \textbf{102} (1991), 95-115.
%
\bibitem{debouard.fukuizumi} A. de Bouard, F. Fukuizumi, \emph{Stability of standing waves for nonlinear Schr\"{o}dinger equations with inhomogeneous nonlinearities}, Ann. Henri Poincar\'{e} \textbf{6} (2005), 1157-1177.
%
%
\bibitem{deng.jin} Y. Deng, L. Jin, S. Peng, \emph{Solutions of Schr\"{o}dinger equations with inverse square potential and critical nonlinearity}, J. Differential Equations \textbf{253} (2012), 1376-1398.
%
%
\bibitem{dupaigne} L. Dupaigne, \emph{A nonlinear elliptic PDE with the inverse square potential}, J. Anal. Math. \textbf{86} (2002), 359-398.
%
\bibitem{ef} T. Ekholm, R. L. Frank, \emph{On Lieb-Thirring inequalities for Schr\"{o}dinger operators with virtual level}, Comm. Math. Phys.  \textbf{264} (2006), 725-740.
%
\bibitem{felli.ferrero} V. Felli, A. Ferrero, \emph{On semilinear elliptic equations with borderline Hardy potentials}, J. Anal. Math. \textbf{123} (2014), 303-340.
%
%
%
\bibitem{filippas.tertikas} S. Filippas, A. Tertikas,  \emph{Optimizing improved Hardy inequalities}, J. Funct. Anal. \textbf{192} (2002), 186-233.
%
\bibitem{genoud} F. Genoud, \emph{Existence and orbital stability of standing waves for some nonlinear Schr\"{o}dinger equations, perturbation of a model case}, J. Differential Equations \textbf{246} (2009), 1921-1943.
%
\bibitem{genoud.stuart} F. Genoud, C.A. Stuart, \emph{Schr\"{o}dinger equations with a spatially decaying nonlinearity: existence and stability of standing waves}, Discrete Contin. Dyn. Syst. \textbf{21} (2008), 137-186.
%
\bibitem{guo.chen} B. Guo, J. Chen, \emph{Orbital stability of standing wave solution for a quasilinear Schr\"{o}dinger equation}, Quart. Appl. Math. \textbf{67} (2009), 781-791.
%
%
\bibitem{jeanjean} L. Jeanjean, S. Le Coz, \emph{An existence and stability result for standing waves of nonlinear Schr\"{o}dinger equations}, Adv. Differential Equations \textbf{11} (2006), 813-840.
%
\bibitem{rose.weinstein} H. A. Rose, M. I. Weinstein, \emph{On the bound states of the nonlinear Schr\"{o}dinger equation with a linear potential}, Phys. D \textbf{30} (1988), 207-218.
%
\bibitem{shibata} M. Shibata, \emph{Stable standing waves of nonlinear Schr\"{o}dinger equations with a general nonlinear term}, Manuscripta Math. \textbf{143} (2014), 221-237.
%
%
\bibitem{smets} D. Smets, \emph{Nonlinear Schr\"{o}dinger equations with Hardy potential and critical nonlinearities}, Trans. Amer. Math. Soc. \textbf{357} (2005), 2909-2938.
%
%
\bibitem{su.wang.willem} J. Su, Z.-Q. Wang, M. Willem, \emph{Nonlinear Schr\"{o}dinger equations with unbounded and decaying radial potentials}, Commun. Contemp. Math. \textbf{9} (2007), 571-583.
%
\bibitem{su.wang.willem.jde} J. Su, Z.-Q. Wang, M. Willem, \emph{Weighted Sobolev embedding with unbounded and decaying radial potentials}, J. Differential Equations \textbf{238} (2007), 201-219.
%
\bibitem{terracini} S. Terracini, \emph{On positive entire solutions to a class of equations with a singular coefficient and critical exponent}, Adv. Differential Equations \textbf{1} (1996), 241-264.
%
%
\bibitem{trach.zogr} G. P. Trachanas, N. B. Zographopoulos, \emph{A strongly singular parabolic problem on an unbounded domain}, Commun. Pure Appl. Anal. \textbf{13} (2014), 789-809.
%
%
\bibitem{vzog1} J. L. V\'azquez, N. B. Zographopoulos, \emph{Functional
aspects of the Hardy inequality: appearance of a hidden energy}, J. Evol. Equ. \textbf{12} (2012), 713-739.
%
\bibitem{vzog2} J. L. V\'azquez, N. B. Zographopoulos,
\emph{Hardy type inequalities and hidden energies},
Discrete Contin. Dyn. Syst. 33 (2013), 5457-5491.
%
\bibitem{vz00} J. L. V\'azquez, E. Zuazua, \emph{The Hardy
inequality and the asymptotic behaviour of the heat equation with
an inverse-square potential},  J. Funct. Anal. \textbf{173} (2000), 103-153.
%
\bibitem{willem} M. Willem, \emph{Minimax theorems}, Progress in Nonlinear Differential Equations and their Applications \textbf{24}, 1996.
%
%
\bibitem{zhang.only} Q. S. Zhang, \emph{Global bounds of Schr\"{o}dinger heat kernels with negative potentials}, J. Funct. Anal. \textbf{182} (2001), 344-370.
%
%
\bibitem{zhang.zheng} J. Zhang, J. Zheng, \emph{Scattering theory for nonlinear Schr\"{o}dinger equations with inverse-square potential}, J. Funct. Anal. \textbf{267} (2014), 2907-2932.
%
%
\bibitem{zogr.jfa} N. B. Zographopoulos, \emph{Existence of extremal functions for a Hardy-Sobolev inequality}, J. Funct. Anal. \textbf{259} (2010), 308-314.
%
\end{thebibliography}

}

\vskip 0.5cm

{\sc Addresses:}

{\sc Georgios P. Trachanas}\newline
Department of Mathematics, National Technical University of Athens, Zografou Campus, 15780, Athens, Greece, \newline
e-mail: gtrach@gmail.com, gtrach@math.ntua.gr

\medskip

{\sc Nikolaos B. Zographopoulos}\newline
University of Military Education, Hellenic Army Academy, Department of Mathematics \& Engineering Sciences, Vari -
16673, Athens Greece,
\newline e-mail: nzograp@gmail.com, zographopoulosn@sse.gr

\

{\bf Keywords:} \ {Orbital Stability, Standing wave, Hardy Inequality.} \vspace{0.1cm}

\medskip

{\bf Mathematics Subject Classification (2010):} \  {35A15, 35B35, 35Q41, 35Q55.}
\end{document}